# On the Total Edge Product Cordial Labeling of Some Corona Graphs

James Cyrile G. Valdehueza, Ariel C. Pedrano

*Department of Mathematics and Statistics, College of Arts and Sciences*
*University of Southeastern Philippines, Davao City, Philippines*

***Abstract -*** *An edge labeling function $f^*: E(G) \to \{0,1\}$ induces a vertex labeling function $f: V(G) \to \{0,1\}$ defined as $f(v) = \prod\{f^*(uv)/uv \in E(G)\}$. The function $f^*$ is called a total edge product cordial labeling of $G$ if $\left|\left(v_f(0) + e_f(0)\right) - \left(v_f(1) + e_f(1)\right)\right| \leq 1$. In this paper, we determine the total edge product cordial labeling of the corona graphs $P_n \circ P_m$ and $P_n \circ C_m$.*

**Keywords** — *Total Edge Product Cordial, Corona graphs, Graph Labeling.*

## I. INTRODUCTION

Graph Labeling has been one of the most interesting topics in the field of graph theory. It is the process of assigning integers in the vertices, some are in the edges while others are both to a certain graph $G$. There are different types of labeling, and one of it is Cordial labeling which was introduced by Cahit [1] in the year 1987. This labeling is said to be a weaker version of graceful and harmonious labeling where we use either 0 or 1 to label the edges and vertices of a certain graph. Moreover, cordial labeling comes in many ways depending upon the kind of conditions being considered.

In 2004, Sundaram et al [2] introduced the notion of product cordial labeling in which the absolute difference in cordial labeling is replaced by the product of the vertex labels. In 2006, Sundaram et al [3] discussed total product cordial labeling of graphs obtained as the join of two graphs. The total product cordial labeling is defined as a binary labeling $f: (V \cup E) \to \{0,1\}$ such that $f(xy) = f(x)f(y)$ where $x, y \in V(G), xy \in E(G)$, and the total number of 0 and 1 are balanced. That is if $v_f(i)$ and $e_f(i)$ denote the set of vertices and edges which are labeled as $i$, $for\ i = 0,1\ repectively$, then $\left|\left(\left|v_f(0)\right| + \left|e_f(0)\right|\right) - \left(\left|v_f(1)\right| + \left|e_f(1)\right|\right)\right| \leq 1$.

In 2012, Vaidya and Barasara [4] introduced the edge analogue of product cordial labeling as a variant of total product cordial labeling and called it total edge product cordial labeling. For a graph $G$, an edge labeling function $f^*: E(G) \to \{0,1\}$ induces a vertex labeling function $f: V(G) \to \{0,1\}$ defined as $f(v) = \prod\{f^*(uv)/uv \in E(G)\}$. The function $f^*$ is called a total edge product cordial labeling of $G$ if $\left|\left(v_f(0) + e_f(0)\right) - \left(v_f(1) + e_f(1)\right)\right| \leq 1$. A graph is called total edge product cordial if it admits total edge product cordial labeling. Various results were already published showing the different graphs that admit total edge product cordial labeling. In the same paper of Vaidya and Barasara [5], it has been revealed that some cycle graphs admit this kind of labeling. Similarly, graphs like wheel, gear, complete, complete bipartite, fan, and double fan admit total edge product cordial labeling under some restrictions.

## II. BASIC CONCEPTS

**Theorem 2.1.** [6] The graph with degree sequences (1,1), (2,2,2,2) or (3,2,2,1) are not total edge product cordial graphs.

**Theorem 2.2.** [6] The wheel $W_m$ is a total edge product cordial graph.

**Theorem 2.3.** [6] The fan $F_m$ is a total edge product cordial graph.

## III. RESULTS AND DISCUSSIONS

This chapter presents some results of Total Edge Product Cordial Labeling of some corona graphs.

**Theorem 3.1.** The corona graph $P_n \circ P_m$ is a total edge product cordial graph for all $n$ and $m$ except when

$n = m = 1$.

*Proof.* Let $V(P_n) = \{u_1, u_2, \ldots, u_n\}$ be the vertex set of $P_n$ and $V(P_m^i) = \{v_1^i, v_2^i, \ldots, v_m^i\}$ be the vertex set of the ith copy of $P_m$. Consider the following cases.

**Case 1: $n = 1$, and $m \geq 2$**

*Subcase 1.1 : $n = 1$ and $m = 1$*

If $n = m = 1$, then $P_1 \circ P_1$ has degree sequence of $(1,1)$. Therefore, by Theorem 2.1, $P_1 \circ P_1$ is not a total edge product cordial graph.

*Subcase 1.2: $n = 1$ and $m \geq 2$*

If $n = 1$ and $m \geq 2$, then $P_1 \circ P_m \cong F_m$. Thus by Theorem 2.3, $P_1 \circ P_m$ is a total edge product cordial graph.

Thus, the corona graph $P_n \circ P_m$ is a total edge product cordial for $n = 1$ and $m \geq 1$, except when $n = m = 1$.

**Case 2: $n$ is even and $m \geq 1$.**

Define the function $f^*: E(P_n \circ P_m) \to \{0,1\}$ as follows:

$$f_1^*(u_i u_{i+1}) = \begin{cases} 0, & 1 \leq i \leq \frac{n}{2} - 1 \\ 1, & \frac{n}{2} \leq i \leq n - 1 \end{cases}$$

$$f_2^*(u_i v_j^i) = \begin{cases} 0, & 1 \leq i \leq \frac{n}{2}\ ; & 1 \leq j \leq m \\ 1, & \frac{n}{2} + 1 \leq i \leq n\ ; & 1 \leq j \leq m \end{cases}$$

$$f_3^*(v_j^i v_{j+1}^i) = \begin{cases} 0, & 1 \leq i \leq \frac{n}{2}\ ; & 1 \leq j \leq m - 1 \\ 1, & \frac{n}{2} + 1 \leq i \leq n\ ; & 1 \leq j \leq m - 1 \end{cases}$$

From the function defined above, the number of edges labeled with 0 and 1 is,

$$e_f(0) = nm - 1 \text{ and } e_f(1) = nm$$

And the number of vertices labeled 0 and 1 is,

$$v_f(0) = \frac{n+nm}{2} \text{ and } v_f(1) = \frac{n+nm}{2}$$

Hence, we have,

$$\left|\left(v_f(0) + e_f(0)\right) - \left(v_f(1) + e_f(1)\right)\right| \leq 1$$

Thus, the corona graph $P_n \circ P_m$ is a total edge product cordial graph if $n$ is even and $m \geq 1$.

**Case 3: $n$ is odd, $n \geq 3$ and $m \geq 1$**

*Subcase 3.1: $n$ is odd, $n \geq 3$ and $m$ is even, $m \geq 1$*

Define the function $f^*: E(P_n \circ P_m) \to \{0,1\}$ as follows:

$$f_1^*(u_i u_{i+1}) = \begin{cases} 0, & 1 \leq i \leq \frac{n-1}{2} \\ 1, & \frac{n+1}{2} \leq i \leq n - 1 \end{cases}$$

$$f_2^*(u_i v_j^i) = \begin{cases} 0, & 1 \leq i \leq \frac{n-1}{2} \quad ; \quad 1 \leq j \leq m \text{ or} \\ & i = \frac{n+1}{2} \quad ; \quad 1 \leq j \leq \frac{m}{2} \\ 1, & \frac{n+3}{2} \leq i \leq n \quad ; \quad 1 \leq j \leq m \text{ or} \\ & i = \frac{(n+1)}{2} \quad ; \quad \frac{m}{2} + 1 \leq j \leq m \end{cases}$$

$$f_3^*(v_j^i v_{j+1}^i) = \begin{cases} 0, & 1 \leq i \leq \frac{n-1}{2} \quad ; \quad 1 \leq j \leq m-1 \text{ or} \\ & i = \frac{n+1}{2} \quad ; \quad 1 \leq j \leq \frac{m}{2} - 1 \\ 1, & \frac{n+3}{2} \leq i \leq n \quad ; \quad 1 \leq j \leq m-1 \text{ or} \\ & i = \frac{(n+1)}{2} \quad ; \quad \frac{m}{2} \leq j \leq m-1 \end{cases}$$

In view of the above labeling, we have, $e_f(0) = nm - 1$ and $e_f(1) = nm$ and the number of vertices labeled 0 and 1 is, $v_f(0) = \frac{nm+n+1}{2}$ and $v_f(1) = \frac{n+nm-1}{2}$, respectively.

Hence, we have,

$$\left|\left(v_f(0) + e_f(0)\right) - \left(v_f(1) + e_f(1)\right)\right| \leq 1$$

Thus, the corona graph $P_n \circ P_m$ is a total edge product cordial graph if $n$ is odd, $n \geq 3$ and $m$ is even, $m \geq 1$.

*Subcase 3.2*: $n$ is odd, $n \geq 3$ and $m$ is odd, $m \geq 1$.

Define the function $f^*: E(P_n \circ P_m) \to \{0,1\}$ as follows:

$$f_1^*(u_i u_{i+1}) = \begin{cases} 0, & 1 \leq i \leq \frac{n-1}{2} \\ 1, & \frac{n+1}{2} \leq i \leq n-1 \end{cases}$$

$$f_2^*(u_i v_j^i) = \begin{cases} 0, & 1 \leq i \leq \frac{n-1}{2} \quad ; \quad 1 \leq j \leq m \text{ or} \\ & i = \frac{n+1}{2} \quad ; \quad 1 \leq j \leq \frac{m-1}{2} \\ 1, & \frac{n+3}{2} \leq i \leq n \quad ; \quad 1 \leq j \leq m \text{ or} \\ & i = \frac{(n+1)}{2} \quad ; \quad \frac{m+1}{2} \leq j \leq m \end{cases}$$

$$f_3^*(v_j^i v_{j+1}^i) = \begin{cases} 0, & 1 \leq i \leq \frac{n-1}{2} \quad ; \quad 1 \leq j \leq m-1 \text{ or} \\ & i = \frac{n+1}{2} \quad ; \quad 1 \leq j \leq \frac{m-1}{2} \\ 1, & \frac{n+3}{2} \leq i \leq n \quad ; \quad 1 \leq j \leq m-1 \text{ or} \end{cases}$$

$$i = \frac{(n+1)}{2} \quad ; \quad \frac{m+1}{2} \leq j \leq m-1$$

From the edge labeling function defined above, we have,

$$e_f(0) = nm - 1 \text{ and } e_f(1) = nm$$

And the number of vertices labeled 0 and 1 is,

$$v_f(0) = \frac{nm+n+2}{2} \text{ and } v_f(1) = \frac{n+nm-2}{2}$$

Hence, we have,

$$\left| \left( v_f(0) + e_f(0) \right) - \left( v_f(1) + e_f(1) \right) \right| \leq 1$$

Thus, the corona graph $P_n \circ P_m$ is a total edge product cordial graph if $n$ is odd, $n \geq 3$ and $m$ is odd, $m \geq 1$.

Considering all the cases above, we can say that the corona graph $P_n \circ P_m$ admits total edge product cordial labeling for all $n$ and $m$ except when $n = m = 1$.

**Theorem 3.2.** The corona graph $P_n \circ C_m$ is a total edge product cordial graph for $n \geq 1$ and $m \geq 3$.

*Proof.* Let $V(P_n) = \{u_1, u_2, \ldots, u_n\}$ be the vertex set of $P_n$ and $V(C_m^i) = \{v_1^i, v_2^i, \ldots, v_m^i\}$ be the vertex set of the ith copy of $P_m$. Consider the following cases.

**Case 1: $n = 1$ and $m \geq 3$.** If $n = 1$ and $m \geq 3$, then $P_n \circ C_m \cong W_m$. Hence, by Theorem 2.2, $P_n \circ C_m$ is a total edge product cordial graph.

**Case 2: $n$ is even and $m \geq 3$**

Define the function $f^*: E(P_n \circ C_m) \to \{0,1\}$ as follows:

$$f_1^*(u_i u_{i+1}) = \begin{cases} 0, & 1 \leq i \leq \frac{n}{2} - 1 \\ 1, & \frac{n}{2} \leq i \leq n - 1 \end{cases}$$

$$f_2^*(u_i v_j^i) = \begin{cases} 0, & 1 \leq i \leq \frac{n}{2} \quad ; \quad 1 \leq j \leq m \\ 1, & \frac{n}{2} + 1 \leq i \leq n \quad ; \quad 1 \leq j \leq m \end{cases}$$

$$f_3^*(v_j^i v_{j+1}^i) = \begin{cases} 0, & 1 \leq i \leq \frac{n}{2} \quad ; \quad 1 \leq j \leq m-1 \\ 1, & \frac{n}{2} + 1 \leq i \leq n \quad ; \quad 1 \leq j \leq m-1 \end{cases}$$

$$f_4^*(v_1^i v_m^i) = \begin{cases} 0, & 1 \leq i \leq \frac{n}{2} \\ 1, & \frac{n+1}{2} \leq i \leq n \end{cases}$$

In view of the above labeling function, the total number of edges labeled with 0 and 1 is,

$$e_f(0) = \frac{2nm+n-2}{2} \text{ and } e_f(1) = \frac{2nm+n}{2}$$

And the total number of vertices labeled with 0 and 1 is,

$$v_f(0) = \frac{n+nm}{2} \text{ and } v_f(1) = \frac{n+nm}{2}$$

Hence, we have,

$$\left|\left(v_f(0)+e_f(0)\right)-\left(v_f(1)+e_f(1)\right)\right|\leq 1$$

Thus, the corona graph $P_n°C_m$ is a total edge product cordial graph if $n$ is even and $m\geq 3$.

**Case 3: $n$ is odd, $n\geq 3$ and $m\geq 3$**

*Subcase 3.1: $n$ is odd, $n\geq 3$ and $m$ is even, $m\geq 3$*

Define the function $f^*: E(P_n°C_m)\to\{0,1\}$ as follows:

$$f_1^*(u_iu_{i+1}) = \begin{cases} 0, & 1\leq i\leq \frac{n-1}{2} \\ 1, & \frac{n+1}{2}\leq i\leq n-1 \end{cases}$$

$$f_2^*(u_iv_j^i) = \begin{cases} 0, & 1\leq i\leq \frac{n-1}{2} \quad ; \quad 1\leq j\leq m \text{ or} \\ & i=\frac{n+1}{2} \quad ; \quad 1\leq j\leq \frac{m}{2} \\ 1, & \frac{n+3}{2}\leq i\leq n \quad ; \quad 1\leq j\leq m \text{ or} \\ & i=\frac{(n+1)}{2} \quad ; \quad \frac{m}{2}+1\leq j\leq m \end{cases}$$

$$f_3^*(v_j^iv_{j+1}^i) = \begin{cases} 0, & 1\leq i\leq \frac{n-1}{2} \quad ; \quad 1\leq j\leq m-1 \text{ or} \\ & i=\frac{n+1}{2} \quad ; \quad 1\leq j\leq \frac{m}{2}-1 \\ 1, & \frac{n+3}{2}\leq i\leq n \quad ; \quad 1\leq j\leq m-1 \text{ or} \\ & i=\frac{(n+1)}{2} \quad ; \quad \frac{m}{2}\leq j\leq m-1 \end{cases}$$

$$f_4^*(v_1^iv_m^i) = \begin{cases} 0, & 1\leq i\leq \frac{n-1}{2} \\ 1, & \frac{n+1}{2}\leq i\leq n \end{cases}$$

Now, the total number of edges labeled 0 and 1 is,

$$e_f(0)=\frac{2nm+n-3}{2} \text{ and } e_f(1)=\frac{2nm+n+1}{2}$$

And the total number of vertices labeled with 0 and 1 is,

$$v_f(0)=\frac{n+nm+1}{2} \text{ and } v_f(1)=\frac{n+nm-1}{2}$$

Hence, we have,

$$\left|\left(v_f(0)+e_f(0)\right)-\left(v_f(1)+e_f(1)\right)\right|\leq 1$$

Thus, the corona graph $P_n°C_m$ is a total edge product cordial graph if $n$ is odd, $n\geq 3$ and $m$ is even, $m\geq 3$.

*Subcase 3.2: $n$ is odd, $n\geq 3$ and $m$ is odd, $m\geq 3$.*

Define the function $f^*: E(P_n°C_m)\to\{0,1\}$ as follows:

$$f_1^*(u_i u_{i+1}) = \begin{cases} 0, & 1 \leq i \leq \frac{n-1}{2} \\ 1, & \frac{n+1}{2} \leq i \leq n-1 \end{cases}$$

$$f_2^*(u_i v_j^i) = \begin{cases} 0, & 1 \leq i \leq \frac{n-1}{2} \quad ; \quad 1 \leq j \leq m \text{ or} \\ & i = \frac{n+1}{2} \quad ; \quad 1 \leq j \leq \frac{m-1}{2} \\ 1, & \frac{n+3}{2} \leq i \leq n \quad ; \quad 1 \leq j \leq m \text{ or} \\ & i = \frac{(n+1)}{2} \quad ; \quad \frac{m+1}{2} \leq j \leq m \end{cases}$$

$$f_3^*(v_j^i v_{j+1}^i) = \begin{cases} 0, & 1 \leq i \leq \frac{n-1}{2} \quad ; \quad 1 \leq j \leq m-1 \text{ or} \\ & i = \frac{n+1}{2} \quad ; \quad 1 \leq j \leq \frac{m-1}{2} \\ 1, & \frac{n+3}{2} \leq i \leq n \quad ; \quad 1 \leq j \leq m-1 \text{ or} \\ & i = \frac{(n+1)}{2} \quad ; \quad \frac{m+1}{2} \leq j \leq m-1 \end{cases}$$

$$f_4^*(v_1^i v_m^i) = \begin{cases} 0, & 1 \leq i \leq \frac{n-1}{2} \\ 1, & \frac{n+1}{2} \leq i \leq n \end{cases}$$

The total number of edges labeled with 0 and 1 is,

$$e_f(0) = \frac{2nm+n-3}{2} \text{ and } e_f(1) = \frac{2nm+n+1}{2}$$

And the total number of vertices labeled with 0 and 1 is,

$$v_f(0) = \frac{n+nm}{2} \text{ and } v_f(1) = \frac{n+nm}{2}$$

Hence, we have,

$$\left| \left( v_f(0) + e_f(0) \right) - \left( v_f(1) + e_f(1) \right) \right| \leq 1$$

Thus, the corona graph $P_n \circ C_m$ is a total edge product cordial graph if $n$ is odd, $n \geq 3$ and $m$ is odd, $m \geq 3$.

Considering all the cases above, we have shown that the corona graph $P_n \circ C_m$ is a total edge product cordial graph for $n \geq 1$ and $m \geq 3$.

## IV. CONCLUSIONS

In this paper we have identified that the graphs $P_n \circ P_m$ and $P_n \circ C_m$ are total edge product cordial graphs adhering the following conditions:

- $P_n \circ P_m$ admits total edge product cordial labeling for all $n$ and $m$ except when $n = m = 1$.
- $P_n \circ C_m$ admits total edge product cordial labeling for all $n \geq 1$ and $m \geq 3$.